\documentclass{amsart}
\usepackage{ifthen, citesort}
\usepackage{amssymb}
\usepackage{epsfig}
\usepackage{enumerate}

\usepackage{graphicx}
\usepackage{psfrag}

\DeclareMathOperator{\tr}{tr}

\DeclareMathOperator{\graph}{graph}
\def\R{\mathbb{R}}

\def\dt{\frac{d}{dt}}

\def\theta{\vartheta}
\def\phi{\varphi}
\def\epsilon{\varepsilon}
\newcommand{\A}[1]{\ifthenelse{#1 = 2}{\lvert A\rvert^{#1}}{\tr A^{#1}}}
\newcommand{\Akl}[2]{\ifthenelse{#1 = 2}%
{\left(\lvert A\rvert^{#1}\right)^{#2}}%
{\left(\tr A^{#1}\right)^{#2}}}

\newtheorem{theorem}{Theorem}[section]
\newtheorem{lemma}[theorem]{Lemma}

\newtheorem{corollary}[theorem]{Corollary}

\theoremstyle{definition}
\newtheorem{definition}[theorem]{Definition}

\theoremstyle{remark}
\newtheorem{remark}[theorem]{Remark}
\newtheorem{problems}[theorem]{Open Problems}

\numberwithin{equation}{section}
\begin{document}

\title{Self-similarly expanding networks to \mbox{curve shortening flow}}

\author{Oliver C. Schn\"urer}
\address{Freie Universit\"at Berlin, Arnimallee 6, 14195 Berlin, Germany}
\curraddr{}
\def\rest{@math.fu-berlin.de}
\email{Oliver.Schnuerer\rest}
\thanks{The authors were supported by SFB 647/B3 ``Raum -- Zeit --
  Materie'': Singularity structure, long-time behavior and dynamics
  of non-linear evolution equations}
\author{Felix Schulze}
\address{Freie Universit\"at Berlin, Arnimallee 6, 14195 Berlin, Germany}
\email{Felix.Schulze\rest}


\date{February 2007}

\dedicatory{}


\begin{abstract}
We consider a network in the Euclidean plane that consists of three 
distinct half-lines with common start points. From that network as
initial condition, there exists a network that consists of three 
curves that all start at one point, where they form $120$ degree 
angles, and expands homothetically
under curve shortening flow. We also
prove uniqueness of these networks.
\end{abstract}

\maketitle

\section{Introduction}\label{intro sec}
\noindent We consider networks in $\R^2$ consisting of three curves with common
start points that evolve under curve shortening flow. As curve shortening
flow is the $L^2$-gradient flow for the length functional, it is natural
to assume that each pair of curves encloses an angle of $120^\circ$ at 
their common start points. If such is the case, we say that the curves 
fulfill the balancing condition. In this paper, we investigate curve 
shortening flow of networks for which the balancing condition 
may be violated initially. More precisely, we study three half-lines
with common start points. We prove that there exists a self-similarly
expanding network that approaches the union of the three half-lines for 
small times, evolves under curve shortening flow, and fulfills the 
balancing condition for positive times. Moreover, we show that such
a solution is unique.
\begin{figure}[htb]
\epsfig{file=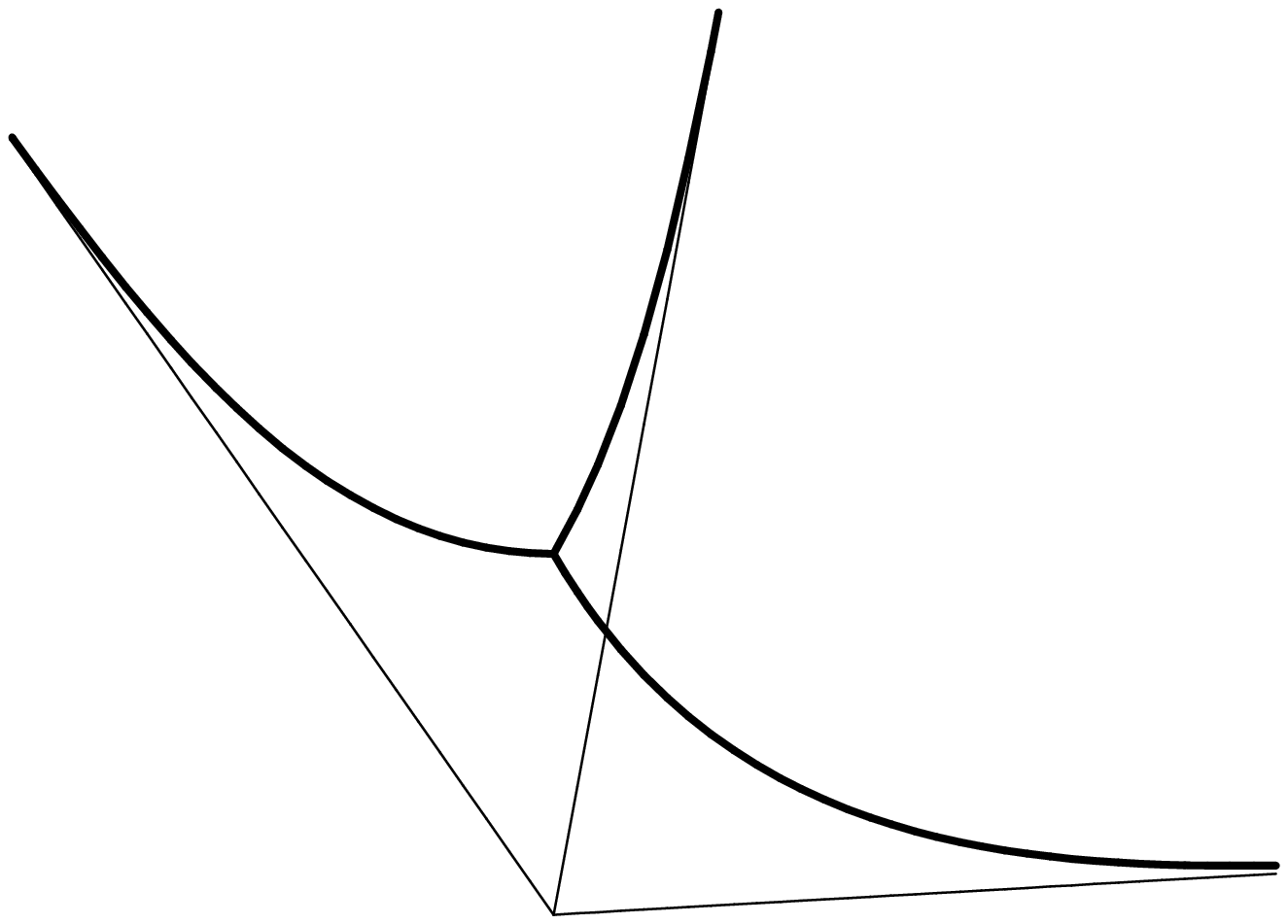, 
  width=0.4\textwidth}
\end{figure}

\noindent In order to formulate our existence and uniqueness result more
precisely, let $l_1$, $l_2$, and $l_3$ be three distinct
half-lines that start in the same point which we assume to be
the origin in $\R^2$. Then there exist directions 
$p_1,\,p_2,\,p_3\in S^1=\left\{x\in\R^2:|x|=1\right\}$ such that
$l_i=\{\lambda p_i:0\le\lambda\}$.
We want to find homothetically expanding solutions to curve
shortening flow for an
initial configuration as described above. Such solutions form
a family of networks $M_t$, $t\ge0$. Each network consists of
three curves $\alpha_i:(\cdot,\,t):[0,\,\infty)\to\R^2$ such that
\begin{itemize}
\item $\alpha_i(r,\,t)$, $i=1,\,2,\,3$, is smooth for $t>0$ 
  and continuous up to $t=0$,
\item for $t>0$, the curves $\alpha_i(\cdot,\,t)$, $i=1,\,2,\,3$, are 
  regular, i.\,e.\ 
  $$\frac{\partial}{\partial r}\alpha_i(r,\,t)\neq0,$$
  uniformly up to $r=0$,
\item the start points $\alpha_i(0,\,t)$, $i=1,\,2,\,3$,
  coincide for all times $t>0$, but may depend on time, 
\item $\alpha_i(\cdot,\,t)$, $i=1,\,2,\,3$, are three embedded curves
  and for all $t$ the three
  curves only meet at the start point, 
\item the tangent vectors at the start point fulfill the
  balancing or $120^\circ$ condition, i.\,e.
  $$\sum\limits_{i=1}^3 T_i=0\text{ for all }t>0,$$
  where
  $$T_i=\frac{\left.\frac{\partial}{\partial r}
  \alpha_i(r,\,t)\right|_{r=0}}
  {\left|\left.\frac{\partial}{\partial r}
  \alpha_i(r,\,t)\right|_{r=0}\right|}$$
  denote the tangent vectors to the curves at $(0,\,t)$, 
\item each curve connects to infinity, i.\,e.
  $$\lim\limits_{r\to\infty}|\alpha_i(r,\,t)|=\infty
  \text{ for all }t\ge0\text{ and }i=1,\,2,\,3,$$
\item the curve $\alpha_i(\cdot,\,t)$, $i=1,\,2,\,3$, $t\ge0$,
  is at infinity asymptotically close to the half-line $l_i$, i.\,e.\
  denoting by $d_{\mathcal H}$ the Hausdorff distance,
  $$\qquad d_{\mathcal H}\left(\alpha_i([0,\infty),t)
  \cap\left(\R^n\setminus B_r(0)\right),\,l_i\cap
  \left(\R^n\setminus B_r(0)\right)\right)
  \to0\quad\text{for }r\to\infty,$$
\item each curve flows for $r>0$ according to curve shortening flow
  (which equals mean curvature or Gau\ss{} curvature flow in
  one space dimension; more precisely, in general, we have to 
  add a non-zero tangential velocity in order to avoid problems
  with the parametrization -- we will later rewrite this flow
  equation equivalently as a graphical flow equation and thereby
  fix a parametrization), i.\,e.
  $$\left(\frac d{dt}\alpha_i\right)^\perp=\Delta\alpha_i,$$
  where $\Delta$ denotes the Laplace-Beltrami operator on the 
  curve,
\item for $t=0$, we get our initial configuration,
  $$\{\alpha_i(r,\,0):r\ge0\}=l_i\text{ for }i=1,\,2,\,3,$$
\item the solution is homothetically expanding, i.\,e.\ for
  $0<t_1<t_2$, there exists $\lambda>1$ such that
  $$\lambda\cdot\bigcup\limits_{i=1}^3\{\alpha_i(r,\,t_1):r\ge0\}
  =\bigcup\limits_{i=1}^3\{\alpha_i(r,\,t_2):r\ge0\},$$
\item $\alpha_i$ is of class $C^0([0,\infty)\times[0,\infty))\cap
  C^\infty([0,\infty)\times(0,\infty))$ for $i=1,\,2,\,3$.
\end{itemize}

\begin{definition}
A family of curves $(\alpha_i)_{i=1,\,2,\,3}$ as described above 
is said to be a network which expands homothetically under
curve shortening flow from three half-lines. 
We say that two such networks are equal if their images
$\bigcup\limits_{i=1}^3\alpha_i([0,\infty),t)$ in $\R^2$ coincide. 
\end{definition}

\noindent In the situation described above, we prove the following
\begin{theorem}\label{straight thm}
For an initial configuration consisting of three distinct half-lines 
starting at the origin, there exists a unique family of networks,
expanding homothetically under curve shortening flow from the
three given half-lines.
\end{theorem}

\noindent The strategy of the proof is as follows. 
For self-similar solutions to curve shortening flow, we can
reduce the parabolic equation to an elliptic equation which
describes self-similar solutions at a fixed positive time,
e.\,g.\ at $t=1/2$. Then we solve this elliptic equation with
prescribed start point $P=\alpha_i(0,1/2)$ such that
$\alpha_i(\cdot,\,1/2)$ connects to infinity and is 
asymptotically close to the half-line $l_i$. This means that 
$M_t=\alpha_i([0,\infty),t)$ is a family of curves that expand
homothetically, is asymptotic to $l_i$ at spatial infinity and
$\alpha_i(0,1/2)=P$. It is not to
be expected that prospective solutions will fulfill the
balancing condition at $P$ for arbitrary $P\in\R^2$. 
In general, these curves $\alpha_i(\cdot,t)$ fulfill all 
the properties required above for a self-similarly expanding
solution but not necessarily the balancing condition
$$\sum\limits_{i=1}^3T_i=0.$$
We will show, however, that for fixed distinct half-lines $l_i$,
the scalar product
$$\left\langle P,\,\sum\limits_{i=1}^3T_i\right\rangle$$
is negative for large values of $|P|$. As the tangent 
vectors $T_i$ depend continuously on $P$, a mapping
degree argument yields that there exists some 
$P_0\in\R^2$ such that 
$$\sum\limits_{i=1}^3T_i=0,$$
i.\,e.\ by choosing $P_0$ as described, we can also ensure that
the balancing condition is fulfilled. 
So we obtain a homothetically expanding network. 

Uniqueness follows from the observation that for a self-similarly
expanding solution that intersects a family of self-similarly expanding
solutions (all solutions are self-similarly expanding under the
same homotheties) that are all asymptotic to a fixed half-line at 
infinity, the angle, at which these curves intersect, is monotone
along the solution mentioned at the beginning. This follows from
applying the Gau\ss{}-Bonnet theorem and the divergence theorem
and uses also the fact that we consider self-similarly expanding
solutions to curve shortening flow. 

The convergence of a smooth, closed and embedded curve in $\R^2$ to a
'round' point was shown by M.\,Gage and R.\,Hamilton \cite{GageHamilton86} and
M.\,Grayson \cite{Grayson89}. An alternative proof was given by
G.\,Huisken in \cite{HuiskenAsian}. Short-time existence for curve
shortening flow of an embedded network of three curves meeting at one
common point, at which the balancing condition is satisfied, was shown by
C.\,Mantegazza, M.\,Novaga, and V.\,Tortorelli
\cite{MantegazzaNetworks}. Provided that no type II singularities
occur and under some further natural assumptions, they show that such a 
network converges to a minimal configuration as $t \rightarrow
\infty$. The problem of the existence of a smooth network, evolving
by curve shortening flow and satisfying the balancing condition with
an initial network not satisfying the balancing condition is raised
in that paper. In the present paper we give a partial answer to this problem.\\
The evolution of embedded networks, satisfying the balancing condition, arises
naturally in two-dimensional multiphase systems, where the evolution of
the interfaces between different phases can be modeled by curve
shortening flow. This comes from the fact that in certain models the
total energy of the interface is proportional to the total length of
the interface, and curve shortening flow is the $L^2$-gradient flow of
the length functional. An example for such a model is given by the growth of
grain boundaries, see e.g.\! \cite{BronsardReitich93}. For further references
and a more detailed discussion, we would like to refer the reader to
the introduction in \cite{MantegazzaNetworks}.

The rest of the paper is organized as follows. We derive the
equation for self-similarly expanding curves in Section \ref{eqn deriv}.
Existence and continuous dependence on the data for a 
self-similarly expanding curve starting at a fixed point at a 
fixed time which is asymptotic to a given half-line is shown in
Section \ref{ex 1 sec}. We use 
this to construct a solution in Section
\ref{ex 3 sec} and prove uniqueness in Section \ref{uniqueness sec}. 
In Appendix \ref{quadrupel sec}, we indicate how to find other
self-similarly expanding solutions and mention some open problems.

We want to thank Bernold Fiedler, Gerhard Huisken, Tom Ilmanen, 
Stefan Liebscher, Roger Moser, Mariel S\'{a}ez and Tatiana Toro for discussions.

\section{Derivation of the Equation}\label{eqn deriv}

\noindent Assume that a curve $\alpha_i(\cdot,\,t)$ is locally represented
as $\graph u(\cdot,\,t)$ over the real line. Then 
$\graph u(\cdot,\,t)$ flows by curve shortening flow, if
\begin{equation}\label{csf graphs}
\dot u=\sqrt{1+u'^2}\left(\frac{u'}{\sqrt{1+u'^2}}\right)'=
\frac{u''}{1+u'^2}\ .
\end{equation}
Our networks are homothetically expanding. We will rotate
our coordinate system such that the start points of 
$\alpha_i(\cdot,\,t)$ lie on the set $\{(0,\,y)\in\R^2:y\ge0\}$.

Assume first that a curve $\alpha_i(\cdot,\,t)$ is represented
as $\graph u(\cdot,\,t)$, $u(\cdot,\,t):\R_{\ge0}\equiv
\{x\in\R:x\ge0\}\to\R$.
For a homothetically expanding solution, the slope of $u$
at $(0,\,u(0,\,t))\in\R^2$ is independent of $t$,
$$\left.\frac{\partial}{\partial x}u(x,\,t)\right|_{x=0}=\text{const}.$$
Our curves expand self-similarly if the graphs at different
times differ by a homothety as follows
\begin{equation}\label{graphs coincide}
\{(x,\,u(x,\,t)):x>0\}=\lambda(t)\cdot
\{(y,\,u(y,\,1/2)):y>0\}
\end{equation}
for some increasing function $\lambda(t):\R_+\to\R_+$. In view
of the geometrically invariant formulation of mean curvature
flow
$$\dt X=-H\nu,$$
where the embedding vector $X$ scales like $\lambda(t)$,
the mean curvature $H$ like $\frac1{\lambda(t)}$, and the
unit normal $\nu$ is scaling invariant, we deduce that
$\dot\lambda\cdot\lambda=1$. As $\lambda(0)=0$, we have
$\lambda(t)=\sqrt{2t}$. So it suffices to describe an
asymptotically expanding solution at some fixed time $t_0$.
We may assume that $t_0=1/2$ and obtain $\dot\lambda(t_0)
=\lambda(t_0)=1$.
We substitute $y=y(x,t)=\frac1{\lambda(t)}x$ and 
differentiate \eqref{graphs coincide} to get
\begin{align*}
\frac{\partial}{\partial t} u(x,\,t)=&\dot\lambda u(y,\,1/2)-\frac{\dot\lambda}{\lambda}
  \frac{\partial}{\partial y}u(y,\,1/2)x,\\
\frac{\partial}{\partial x}u(x,\,t)=&\frac{\partial}{\partial y}u(y,\,1/2),\\
\frac{\partial^2}{\partial
  x^2}u(x,\,t)=&\frac1{\lambda(t)}\frac{\partial^2}{\partial y^2}
  u(y,\,1/2).
\end{align*}
We use these relations in \eqref{csf graphs} and obtain
$$\dot\lambda(t)\left(u(y,\,1/2)-\frac x{\lambda(t)}
\frac{\partial}{\partial y}u(y,\,1/2)\right)=\frac1{\lambda(t)}
\left(\frac{\frac{\partial^2}{\partial y^2}u(y,\,1/2)}
{1+\left(\frac{\partial}{\partial y}u(y,\,1/2)\right)^2}\right),$$
where $x=\lambda(t)y$. At $t=1/2$, we get
$$u-xu'=\frac{u''}{1+(u')^2}.$$
Note that at $t=1/2$, we have $x=y$ and
$u'=\frac{\partial}{ \partial x}u=\frac{\partial}{\partial y}u$.
So we study the initial value problem
\begin{equation}\label{curve ode}
u-xu'=\frac{u''}{1+(u')^2},\quad x>0,
\end{equation}
for given $u(0)>0$ and $u'(0)$. If $u$ is a solution to
\eqref{curve ode}, then so is $-u$. Solutions with $u(0)=0$
are straight half-lines. Thus it suffices to study solutions
with $u(0)>0$.

\section{Existence of one Curve}\label{ex 1 sec}

\noindent In this section, we show that for every point $P\in\R^2$
and every half-line $l_1$, there exists a self-similarly 
expanding curve 
as described before Theorem \ref{straight thm} that starts
at $P$ and is asymptotic to $l_1$ at infinity. 
This curve $\alpha_1$ fulfills all the conditions mentioned
before Theorem \ref{straight thm} that can be fulfilled by a single
curve.

\begin{lemma}\label{sol convex}
As long as a solution to \eqref{curve ode} exists, i.\,e.\ on
a maximal interval $0\le x<x_{max}$, $u$ is a strictly convex
function.
\end{lemma}
\begin{proof}
Standard theory for ordinary differential equations implies that
a unique smooth solution $u$ to \eqref{curve ode} (with initial
conditions for $u(0)>0$ and $u'(0)$) exists on a maximal interval
$0\le x<x_{max}\le\infty$. Our initial conditions with
$u(0)>0$ and the differential equation ensure that $u''(0)>0$.
Suppose that there exists $x_0>0$ such that $u''(x_0)=0$. 
We then have $u(x_0)-x_0u'(x_0)=0$. The point $(x_0,u(x_0))$ 
lies on some straight line 
through the origin. At $x_0$, the slopes of that
straight line and of $u$ coincide. Note that this straight line
solves our differential equation for all $x\in\R$ and is the
only solution $w$ on a maximal existence interval of that
equation with $w(x_0)=u(x_0)$ and $w'(x_0)=u'(x_0)$. As solutions
to \eqref{curve ode} (with given initial conditions) are unique, 
we see that $\graph u$ is contained in a straight line, 
a contradiction.
\end{proof}

\begin{lemma}\label{longtime existence}
The initial value problem \eqref{curve ode} has a solution
for all $x\ge0$.
\end{lemma}
\begin{proof}
If $u'(0)$ is negative, the convexity of $u$ implies that
$|u'|$ is decreasing. Thus $|u'|$ is bounded 
as long as $u'$ is negative. So, for 
bounded values of $x$, $u$ is bounded there. We deduce
that $u$ exists for all $x\ge0$ or $u'$ becomes non-negative.
We may assume the latter.
\par
So we can find $\epsilon>0$ in the maximal existence interval
such that $u'(\epsilon)\ge0$. The convexity of $u$ ensures
that $u'(x)\ge0$ for $x>\epsilon>0$. As $u$ is a strictly
convex function, the differential equation implies that
$u-xu'>0$. So we deduce that
$$u>xu'\ge\epsilon u'\quad\text{for~}x>\epsilon>0.$$
Therefore, according to the maximum principle for ordinary
differential equations, $u$ can grow at most exponentially 
for $x>\epsilon$. Moreover, $u>\epsilon u'$ gives also a bound
on $u'$. Thus $u$ exists for all $x\ge0$.
\end{proof}

\noindent As solutions to \eqref{curve ode} exist for all $x\ge0$,
we will assume in the following that any solution $u$ exists
on an interval containing $[0,\,\infty)$. Note that if 
$u:[0,\infty)\to\R$ is a solution to \eqref{curve ode}, then
$u(-x)$ is a solution to \eqref{curve ode} with 
$u:(-\infty,0]\to\R$. Thus it suffices to study properties
of $u$ on $[0,\infty)$. 

In order to prove that $\graph u$ converges to a straight 
line through the origin for $x\to\infty$, it is useful to 
consider the quantity $u-xu'$ which vanishes precisely 
when $\graph u$ is contained in such a straight line.

\begin{lemma}\label{straight line detection}
For a solution $u$ to \eqref{curve ode}, $u-xu'$ is 
exponentially decaying in $x$.
\end{lemma}
\begin{proof}
As $u$ is strictly convex, the equation implies that
$u-xu'>0$. We use \eqref{curve ode} to derive an evolution 
equation for $u-xu'$ on $\{x\ge\epsilon>0\}$
$$(u-xu')'=-xu''=-x(u-xu')\left(1+(u')^2\right)
\le-\epsilon(u-xu').$$
Therefore $u-xu'$ is exponentially decreasing in $x$.
\end{proof}

\noindent Near $x=0$, the graph of a solution $u$ intersects lines through
the origin of the form $\{(x,\,ax):x\in\R\}$, $a\gg1$. 
By convexity, $u'$ is bounded below, so $u$ will never intersect
lines of the form $\{(x,\,ax):x\in\R\}$ for $a\ll-1$. 
In order to find the line where $\graph u$ becomes 
asymptotic to, we show in the next corollary that $\graph u$
intersects each line through the origin at most once. We will
then prove that for large values of $x$, 
$\graph u$ is asymptotically close to the
line through the origin with slope equal to the infimum 
of slopes of lines through the origin which intersect
$\graph u$.

\begin{lemma}\label{ux/x monoton}
The function $x\mapsto\frac{u(x)}x$ is strictly decreasing.
\end{lemma}
\begin{proof}
We have
$$\left(\frac{u(x)}x\right)'=\frac{u'}x-\frac u{x^2}=
\frac{xu'-u}{x^2}<0.$$  
The Lemma follows.
\end{proof}

\noindent This implies in particular
\begin{corollary}
Let $a\in\R$ and let $u$ be a solution to \eqref{curve ode} with
$u(0)>0$. If $u(x_0)=ax_0$ for some $x_0>0$, then  
$u(x)>ax$ for $0\le x<x_0$ and $u(x)<ax$ for $x>x_0$. Thus
the graphs of $u$ and $x\mapsto ax$ intersect at most once for
$x\ge0$.  
\end{corollary}

\noindent The following lemma implies for some $a>0$ that 
$\graph u|_{\R_{\ge0}}$ is asymptotic to a half-line of the form
$\graph\,(\R_{\ge0}\ni x\mapsto ax)$ at infinity. 
\begin{lemma}\label{exp conv}
There exists $a\in\R$ such that $u(x)-ax\to0$ exponentially
in $C^\infty$ for $x\to\infty$.
\end{lemma}
\begin{proof}
We define $a:=\lim\limits_{x\to\infty}\frac{u(x)}x$. As $u$ is
convex, $u$ decays at most linearly. Thus Lemma \ref{ux/x monoton}
implies that $a\in\R$ is well-defined. 
According to Lemma \ref{straight line detection},
$u(x)-xu'(x)$ converges exponentially to zero for $x\to\infty$.
Therefore we have $\lim\limits_{x\to\infty}u'(x)
=\lim\limits_{x\to\infty}\frac{u(x)}{x}=a$. In view of
Equation \ref{curve ode} we obtain that $u''(x)\to0$ 
exponentially as $x\to\infty$. Integrating yields that
also $u'(x)\to a$ exponentially for $x\to\infty$.
Using Lemma \ref{straight line detection} once again,
we get that $u(x)-ax\to0$ exponentially in $C^0$ as $x\to\infty$.
Differentiating \eqref{curve ode} yields bounds on higher
derivatives of $u$ and interpolation implies exponential
convergence.
\end{proof}

\noindent Solutions to \eqref{curve ode} fulfill useful monotonicity
properties. We will derive these in the following lemmata.

\begin{lemma}\label{height mon}
Let $u$ and $v$ be solutions to \eqref{curve ode} with 
$u(0)>v(0)>0$ and $u'(0)=v'(0)$. Then $u-v$ is positive 
everywhere. This difference is even strictly increasing
in $x$ for $x\ge0$.
\end{lemma} 
\begin{proof}
The differential equation implies that $u''(0)>v''(0)$.
Thus the difference $u-v$ is strictly increasing in $x$
for small $x>0$. Consider a maximal interval 
$(0,\,x_0)$, $0<x_0\le\infty$, where $u'-v'>0$.
We want to exclude that $x_0<\infty$. If this is the
case, we get $u'(x_0)-v'(x_0)=0$ and 
$u''(x_0)-v''(x_0)\le0$. As $u-v$ is strictly increasing
on $(0,\,x_0)$, we have $u(x_0)>v(x_0)$. We obtain
\begin{align*}
0\ge&u''(x_0)-v''(x_0)\\
=&\left(1+(u'(x_0))^2\right)(u(x_0)-x_0u'(x_0))-
\left(1+(v'(x_0))^2\right)(v(x_0)-x_0v'(x_0))\\
=&\left(1+(u'(x_0))^2\right)(u(x_0)-v(x_0))>0,
\end{align*}
a contradiction. Thus $x_0=\infty$ and $u-v$ is strictly
increasing in $x$.
\end{proof}

\begin{lemma}\label{slope mon}
Let $u$ and $v$ be solutions to \eqref{curve ode} with
$u(0)=v(0)>0$ and $u'(0)>v'(0)$. Then $u-v$ is positive
for $x>0$. This difference is even strictly increasing
in $x$ for $x\ge0$. 
\end{lemma}
\begin{proof}
By assumption, $u-v$ is strictly increasing for small $x\ge0$.
Suppose that there is a smallest $x_0>0$ such that 
$u'(x_0)-v'(x_0)=0$. At $x_0$, we have $u(x_0)>v(x_0)$ and
$u''(x_0)\le v''(x_0)$. We deduce there 
$$u-xu'=\frac{u''}{1+(u')^2}\le\frac{v''}{1+(v')^2}=v-xv',$$
so $u(x_0)\le v(x_0)$, a contradiction.
\end{proof}

\noindent Not only for the functions, but also for the asymptotic slopes,
we obtain a strict monotonicity. 

\begin{lemma}\label{different slopes}
Let $u$ and $v$ be solutions to \eqref{curve ode} with
$u(0)\ge v(0)>0$ and $u'(0)\ge v'(0)$, but 
$(u(0),u'(0))\neq(v(0),v'(0))$. Then 
$$\lim\limits_{x\to\infty}\frac{u(x)}x>
\lim\limits_{x\to\infty}\frac{v(x)}x.$$
\end{lemma}
\begin{proof}
Use Lemmata \ref{height mon} and \ref{slope mon} to see that
$u(x)-v(x)$ is positive and increasing in $x$ for $x>0$. As the
exponential convergence in Lemma \ref{exp conv} is uniform for
initial values $(u(0),u'(0))$ in a compact set, the 
limits of $\frac{u(x)}x$ and $\frac{v(x)}x$ differ.
\end{proof}

\noindent For solutions $u$ to \eqref{curve ode} with $u'(0)=0$,
we want to study the slopes $u'$ near infinity. For $u(0)=0$,
we obtain that $u(x)=0$ is a solution. In the following
lemma, we prove that the slopes near infinity become
large for large values of $u(0)$.

\begin{lemma}\label{big slope}
For given $a>0$, there exists $h\gg1$, such that a solution
$u$ to \eqref{curve ode} with $u'(0)=0$ and $u(0)\geq h$ has slope
bigger than $a$ near infinity,
$$\lim\limits_{x\to\infty}u'(x)>a.$$
\end{lemma}
\noindent According to Lemma \ref{ux/x monoton} and the proof of 
Lemma \ref{exp conv}, we also obtain
that $\lim\limits_{x\to\infty}\frac{u(x)}{x}>a$ and
$u(x)>ax$ for $x\ge0$. 
\begin{proof}
If $u'\left(\frac12\right)\ge h$ for some $h\gg1$, 
the lemma is proven. Details for this argument can be found 
at the end of this proof.
Otherwise we have $0\le u'(x)\le h$ and $h\le u(x)$
for all $0\le x\le\frac12$ as $u$ is convex. We estimate there
$$u''=\left(1+(u')^2\right)\cdot(u-xu')\ge 
1\cdot\left(h-\tfrac12 h\right)=\tfrac12h.$$
Thus $$u'\left(\tfrac12\right)\ge\int\limits_0^{1/2}
u''(x)\,dx\ge\tfrac14h.$$
By convexity and Lemma \ref{exp conv},
$\lim\limits_{x\to\infty}u'(x)$ exists and
$\lim\limits_{x\to\infty}u'(x)>u'\left(\frac12\right)\ge\frac14h$.
\end{proof}

\noindent We will later use the following variant of Lemma \ref{big slope}.
\begin{lemma}\label{points inside}
Let $a>0$, $\mu>0$. Then there exists $h\gg1$ such that
every solution $u$ to \eqref{curve ode} with $u(0)\ge h$ and
$\lim\limits_{x\to\infty}u'(x)\le a$ fulfills $u'(0)\le-\mu$.
\end{lemma}
\begin{proof}
As the asymptotic slope is monotone in the height and in 
the initial slope, see Lemma \ref{different slopes}, 
it suffices to show that a solution to
\eqref{curve ode} with $u(0)=h$ and $u'(0)=-\mu$ has an 
asymptotic slope bigger than $a$ at infinity, if $h\gg1$
is chosen sufficiently large. 

According to Lemma \ref{different slopes}, we may assume that
$\mu$ is big; we assume that $\mu>a$. 
As long as $u\ge\frac h2$ on the interval $0\le x\le 1$,
we estimate
$$u''(x)=\left(1+(u'(x))^2\right)(u(x)-xu'(x))
\ge u(x)-xu'(x)\ge\frac h2-|u'(x)|.$$
If $u'(x)>\mu$ for some $x>0$, convexity implies that the
lemma holds. Otherwise, we may assume that $\frac h4\ge\mu$
and use convexity to estimate further
$$u''(x)\ge\frac h2-\mu\ge\frac h4.$$
A solution to $v''(x)=\frac h4$ with initial conditions
$v(0)=h$, $v'(0)=-\mu$ is given by
$v(x)=h-\mu x+\frac h8x^2$. We have $u(x)\ge v(x)$, for
$0\le x\le 1$ provided that $u(x)\ge\frac h2$.
The minimum of $v$ is attained at $x=\frac{4\mu}h$ and
equals $h-\frac{2\mu^2}h$ which is bigger
than $\frac h2$ if $h>2\mu$. We will assume this in
the following. This justifies the assumption that 
$u\ge\frac h2$ on the interval $0\le x\le 1$ for large 
values of $h$. We obtain that $u(x)\ge v(x)$ for $0\le x\le 1$. 
We deduce especially that $u(1)\ge v(1)\ge \frac98h-\mu$.
Assume that $h\ge8\mu$. Then $u(1)-u(0)\ge8\mu$. Therefore,
we find $x\in[0,\,1]$ with $u'(x)\ge8\mu$. Now the Lemma
follows as $u$ is convex.
\end{proof}

\noindent In the following Lemma, we show that for a prescribed half-line
$\{(x,\,ax):x\in\R\}$, $a>0$, there exists an initial height $h$
such that the corresponding solution to \eqref{curve ode} 
converges to that half-line near infinity.

\begin{lemma}\label{exist fixed h}
Let $a>0$. Then there exists an initial height $u(0)=h$ such
that the corresponding solution to \eqref{curve ode} with
$u'(0)=0$ fulfills 
$$\lim\limits_{x\to\infty}u(x)-ax=0.$$
\end{lemma}
\begin{proof}
Apply graphical mean curvature flow to the initial value
$u_0(x)=a\cdot|x|$ and use the results of 
\cite{StavrouSelfSim,EckerHuiskenInvent,EckerHuiskenAnn}.
If $u:\R\times[0,\infty)$ is such a solution, then $u(\cdot,1/2)$
is as claimed in the lemma. 
\end{proof}

We can also prove that $a$ and $h$ depend continuously on each
other.
\begin{lemma}\label{homeo lem}
Let $a$ and $h$ be related via solutions $u$ to \eqref{curve ode}
as in Lemma \ref{exist fixed h}. By reflection of $\graph u$, 
we extend this to a map $\R\ni h\mapsto a\in\R$. Then this map
is a homeomorphism with $a(h)\to\pm\infty$ for
$h\to\pm\infty$. 
\end{lemma}
\begin{proof}
Lemma \ref{big slope} ensures that $a(h)\to\pm\infty$ for
$h\to\pm\infty$. Surjectivity follows from Lemma \ref{exist fixed h}. 
Lemma \ref{different slopes} implies the injectivity. Continuity of
the map follows from the continuous dependence on initial data of solutions 
to ordinary differential equations, applied to $x$ in compact
intervals, and the uniform exponential convergence, 
Lemma \ref{exp conv}. Now standard topology, see e.\,g.\
\cite[Satz 8.12 and Satz 8.24]{Querenburg}, implies that a bijective,
continuous and proper map from $\R^n \to \R^n$ is in fact a
homeomorphism. Thus $h\mapsto a(h)$ is a homeomorphism.
\end{proof}

\noindent We assume now that we are given a half-line $l$ (starting at the
origin) in $\R^2$ and a point $P$. Our aim is to find a 
self-similarly expanding solution to curve shortening flow
with start point at $t=1/2$ equal to $P$ which is asymptotic
to the given half-line $l$. It suffices to consider the elliptic
problem derived in Section \ref{eqn deriv}. We may rotate the
situation such that $P$ lies in the set $\{(0,\,y):y\ge0\}$.
If $P$ lies on the half-line $l$, the claim is obvious. 
Assume thus that $P\neq0$ and that the half-line is given
as $\{(x,\,ax):x\ge0\}$, $a\in\R$.

The idea to find the corresponding curve is as follows. For
$h\in\R$, we consider solutions to \eqref{curve ode} with 
$u(0)=h$ and $u'(0)=0$ for all $x\in\R$. We rotate these 
solutions around the origin such that the half-line
$\graph\,(\R_{\ge0}\ni x\mapsto ax)$ with $a$ as above is mapped
to the given half-line $l$. The resulting curves for different
values of $h$ are visualized in Figure 1.
\begin{figure}[htb]
\epsfig{file=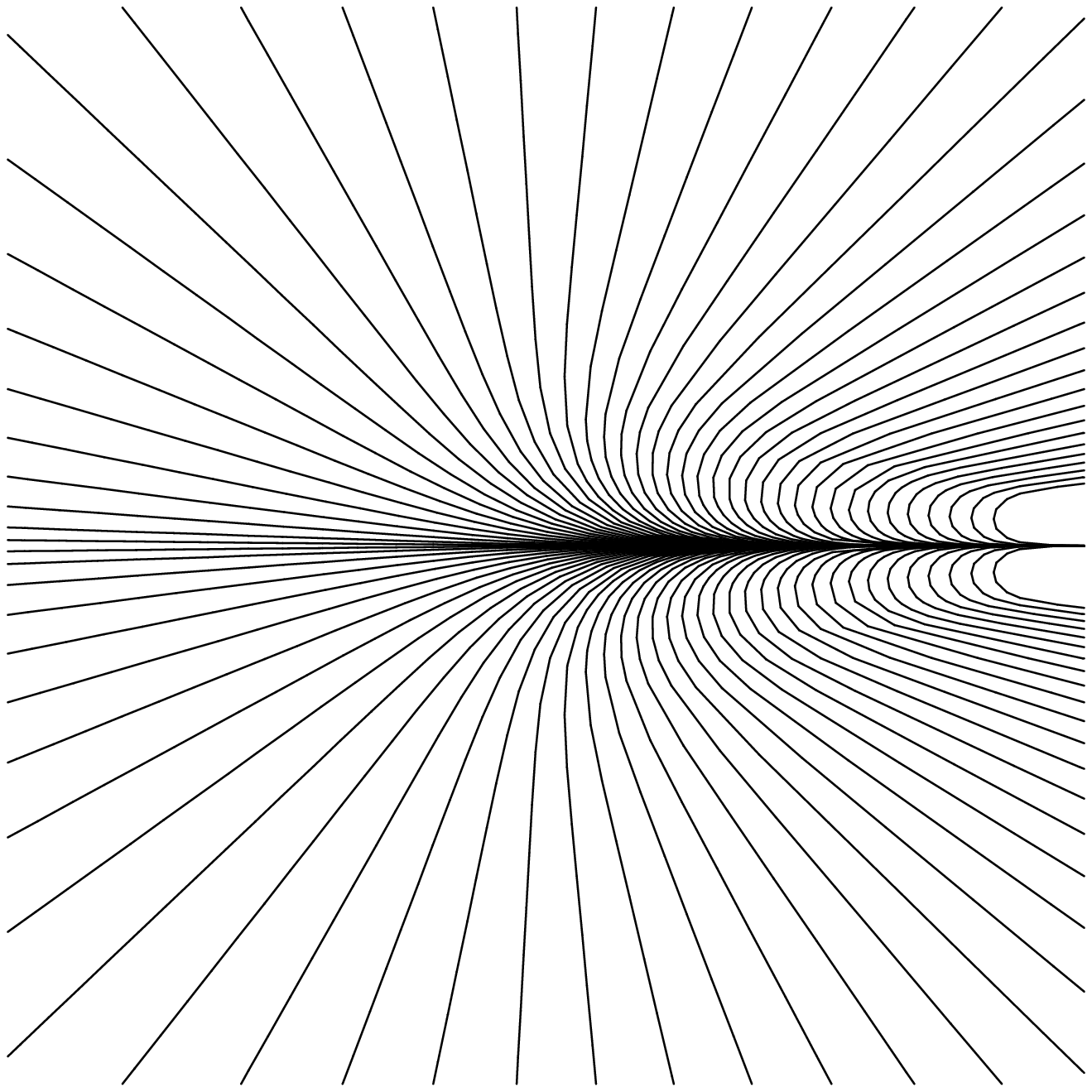, 
  width=0.4\textwidth}
\caption{}
\end{figure}
In the picture, these curves foliate $\R^2$ nicely. So it suffices
to pick the curve that passes through $P$. Note also that for
different values of $P$, the corresponding tangent vectors to
the respective curve seem to vary continuously.

The facts visualized here, are proved as follows. Lemma 
\ref{homeo lem} implies that the rotations depend continuously on
$h$, so we get a family of curves parametrized by $h\in\R$, that 
depend continuously (measured in $C^0_{loc}$ for appropriate
parametrizations) on $h$. For $|h|$ sufficiently large, due 
to Lemma \ref{big slope} and the convexity, the
rotated curves are contained in a cone of small opening angle
around the half-line $l$. We may assume that $P$ has a neighborhood
disjoint to that cone. When we vary $h$, one of the half-lines to
which a curve is asymptotic to at infinity always coincides 
with the half-line $l$ while the other turns by $360^\circ$. This
is only possible if there exists a curve in the family 
considered, that contains $P$ and is asymptotic to $l$. Applying
Lemma \ref{slope mon} and the exponential convergence, Lemma
\ref{exp conv}, similarly to the proof of Lemma \ref{homeo lem},
we see that the curve is also unique. Applying again
\cite[Satz 8.12 and Satz 8.24]{Querenburg} to a map induced from the family 
of curves considered, we see that $(x,h)\in\R^2$ such that 
$(x,u(x))$, where $u$ solves \eqref{curve ode} with 
$u(0)=h$ and $u'(0)=0$, is rotated to $P$, depends 
continuously on $P$. This implies in particular that the tangent
vector to that curve at $P$ depends continuously on $P$. By 
direct inspection, we see that such is also true if $P$ lies
on $l$. We have thus proved
\begin{lemma}
Let $P\in\R^2$ and let $l$ be a half-line starting at the
origin. Then there exists a solution $u$ to \eqref{curve ode}
and a rotation $R$ about the origin such that $P\in
R\,(\graph u)$ and $R\left(\graph u|_{\R_{\ge0}}\right)$ 
is asymptotic to $l$
at infinity. Assume that $P=R\,(x_0,u(x_0))$. Then the tangent
vector $R\Big(\frac{(1,u')}{\sqrt{1+(u')^2}}(x_0)\Big)$ depends 
continuously on $P$.
\end{lemma}

\noindent It also follows that the tangent vector varies continuously
if we rotate the fixed half-line. We do not use this 
observation in the present paper. 

\section{Existence of Three Curves}\label{ex 3 sec}

\noindent In this section, we show that there exists at least one start point 
$P\in\R^2$ for which the problem described in Section \ref{intro sec}
can be solved. 

We have seen in Section \ref{ex 1 sec}, that for any
$P\in\R^2$, there exist three curves that fulfill the
conditions in Section \ref{intro sec}, but do not
necessarily meet at the start point at $120^\circ$ angles.
We also have to show that these curves only have their
start point in common.

The following lemma concerns embeddedness of the curves.
\begin{lemma}
Different self-similarly expanding curves which have
a common start point do not intersect outside this 
start point.
\end{lemma}
\begin{proof}
We may rotate our coordinate system such that the start
point lies on the half-line $\{(0,\,y):y\ge0\}$. 
These curves are always graphs over $\{x\ge0\}$
or over $\{x\le0\}$ or are part of the axis $\{x=0\}$.
Thus it suffices to study two curves that are represented as
graphs over the half-line $\{x\ge0\}$. They start at the
same height and have different initial slopes as they are
not identical. Thus Lemma \ref{slope mon} implies that
these two curves meet only for $x=0$. 
\end{proof}

\noindent In order to show that there exists at least one start point $P$,
where the three curves meet at an angle of $120^\circ$,
we fix the half-lines and vary $P$. For each $P$, we sum
up the tangential directions $(T_i)_i$ at $P$. This gives a 
continuous vector field $V=\sum\limits_{i=1}^3T_i$ on $\R^2$. 
We will show that for points $P$
far from the origin, we always get $\langle P,\, V(P)\rangle<0$.
Thus, by the Brouwer fixed point theorem, there exists some $P_0$ 
such that $V(P_0)=0$. For this start
point, the $120^\circ$ condition is fulfilled. 

The proof of Theorem \ref{straight thm} is thus finished
provided that we can show that $\langle P,V(P)\rangle<0$
for all $P\in\R^2$ such that $|P|$ is sufficiently large. 
This follows from
Lemma \ref{points inside}. It states that for $|P|$ sufficiently
large, a self-similar solution has to start almost in the direction
towards the origin in order to be asymptotic to a half-line that
encloses an angle estimated from below (by a positive constant)
to the half-line starting at the origin and passing through $P$. 
Thus the corresponding tangent vector at $P$ almost equals
$-\frac P{|P|}$.
Thus, having fixed three distinct half-lines starting at the
origin, the half-line starting at the origin and passing through
$P$ encloses an angle estimated from below to at least two of
them. Thus at least two of the three vectors $T_i$ 
at $P$ are
almost equal to $-\frac P{|P|}$. This suffices to ensure that
$\langle P,\,V(P)\rangle=\sum\limits_{i=1}^3\langle P,\,T_i\rangle
<0$. Therefore, we find some $P\in\R^2$ such that we have
$\sum\limits_{i=1}^3T_i=0$ there, i.\,e.\ every two tangent vectors
at $P$ enclose an angle of $120^\circ$. This completes the
existence part of the proof of Theorem \ref{straight thm}.

\section{Uniqueness}\label{uniqueness sec}

\noindent To prove uniqueness of a homothetically expanding network for three
given distinct half-lines, we will first show that the angle, at
which self-similarly expanding curves, asymptotic to a fixed half-line,
intersect another self-similarly expanding curve, is monotone.

Let $l$ be a given-half line in $\mathbb{R}^2$, which we can assume to
coincide with the set $\{(x,0)\in \mathbb{R}^2\ | \ x\geq
0\}$. Furthermore let 
\mbox{$c_h:\R\rightarrow \mathbb{R}^2$} be the self-similarly expanding curves
asymptotic to this half-line, as considered in Section \ref{ex 1
  sec}, parametrized such that
$$\text{dist}(c_h(t), l) \rightarrow 0 \quad \text{as}\ t \rightarrow \infty$$
 We distinguish these curves by the parameter $h$,
where $h$ is such that 
$$\text{dist}\left(\text{Im}(c_h),(0,0)\right) = |h|\  ,$$
and we take $h$ to be positive if $c_h$ lies in the upper half-plane
and negative, if $c_h$ lies in the lower half-plane.\\
Let $b:\mathbb{R}\rightarrow \mathbb{R}^2$ be a further self-similarly
expanding curve, asymptotic to a different half line in $\mathbb{R}^2$,
i.e. there exists an $h\in \mathbb{R}$ and a rotation $R\neq
\text{id}$ such that $b = R\circ c_h$. Note that by Lemma \ref{slope
  mon}, a curve $c_h$ intersects $b$ at most once. Assume now that we
have two curves $c_{h_1},c_{h_2}$ which intersect $b$
in points $p_1$ and $p_2$, respectively. We assume further that
$h_1,h_2$ are such that $p_1$ lies before $p_2$ with respect to the
parametrization of $b$.
Let $\alpha_i$ be the angle
between the tangent vector $T_i$ to $c_{h_i}$ and the tangent vector
$T_b$ at $p_i$, where $i=1,2$. 

\begin{lemma}\label{angle mon} We have
$$ \alpha_2-\alpha_1 = 2 |A|, $$
 where $A$ is the triangle enclosed by these three curves as indicated
 in the picture and $|A|$ denotes its area.
\end{lemma}

\noindent Note that by the exponential convergence of $c_{h_1}(t)$ and
$c_{h_2}(t)$, parametrized as graphs over the $x$-axis, to
$l$ as $t \rightarrow \infty$, the area $|A|$ is finite.

\begin{figure}[htb]
\psfrag{a}{$\alpha_1$}  
\psfrag{b}{$\alpha_2$}  
\psfrag{c}{$c_{h_2}$}   
\psfrag{d}[r]{$c_{h_1}$}
\psfrag{e}{$b$}
\psfrag{f}{$A$}
\epsfig{file=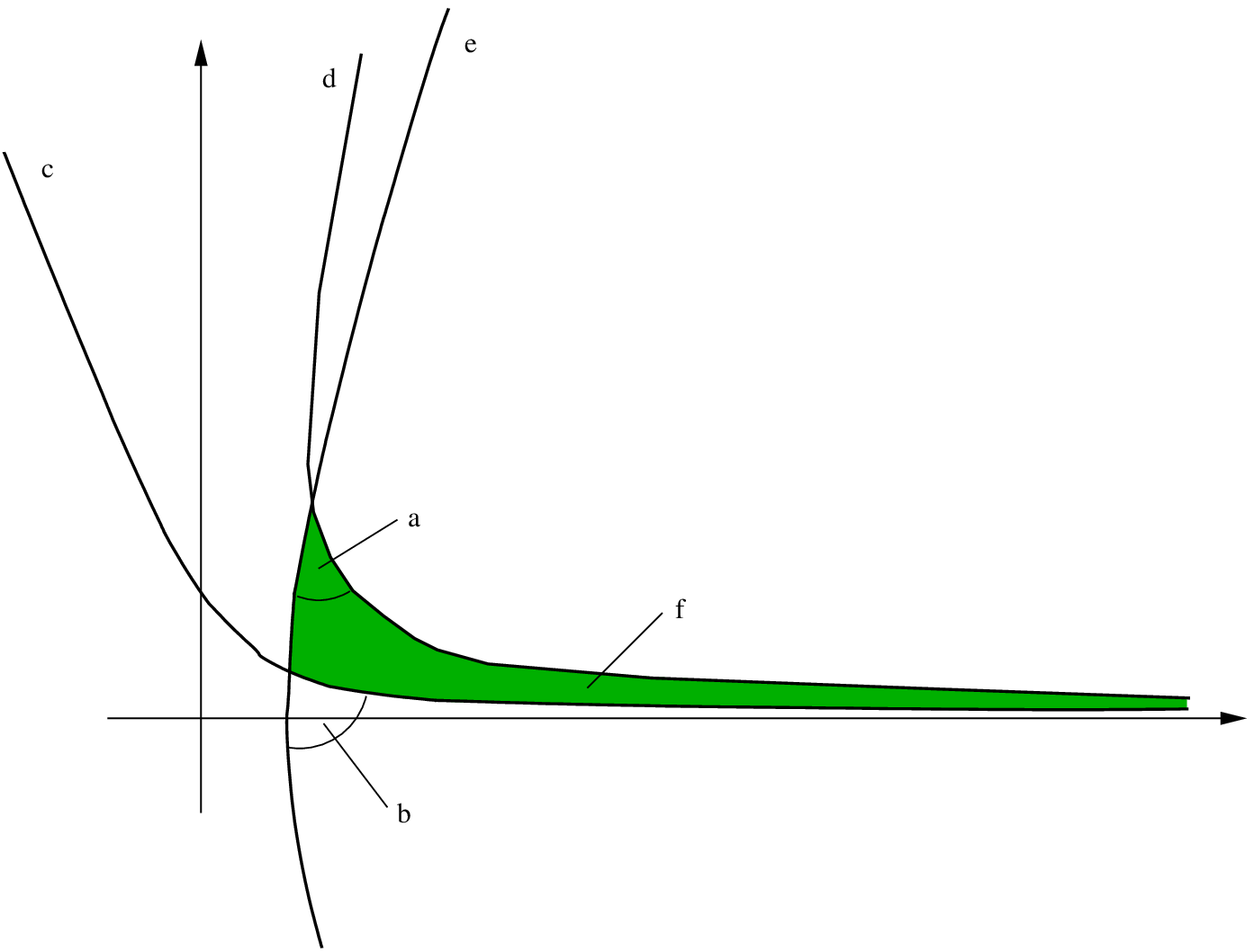, 
  width=0.6\textwidth} 
\end{figure}

\begin{proof} By an analogous derivation as in Section \ref{eqn deriv}
  it is easy to see that equation \eqref{csf graphs} for a
  self-similarly expanding curve can be written in the coordinate free
  form
$$H= - \langle X, \nu \rangle\ ,$$
where the sign of $H$ is chosen such that $-H\nu$ is the (mean)
curvature vector, $\nu$ is a choice of unit normal and $X$ is the
position vectorfield $(x_1, x_2)\mapsto (x_1,x_2)$. By the
divergence theorem we can compute
\begin{align*}
2|A|=&\ \int\limits_A \text{div}(X)\, dx = \int\limits_{\partial A} \langle X,
\nu\rangle\, d\sigma = -\int\limits_{\partial A} \kappa\, d\sigma\\
=&\ -\left(2\pi - \pi - (\pi-\alpha_1) -\alpha_2\right) =
\alpha_2-\alpha_1\ ,
\end{align*}
where $\nu$ is the outward unit normal to $\partial A$ and we have
oriented $\partial A$ counterclockwise. 
\end{proof}

\begin{figure}[htb]
\psfrag{a}{$\eta$}  
\psfrag{b}[l]{$\gamma$}  
\psfrag{c}{$\zeta$}  
\psfrag{d}[c]{$b_3$}
\psfrag{e}{$b_1$}
\psfrag{f}[r]{$b_2$}
\psfrag{g}[l]{$c_1$}
\psfrag{h}[r]{$c_2$}
\psfrag{i}[l]{$c_3$}
\psfrag{j}{$l_1$}
\psfrag{l}[l]{$l_2$}
\psfrag{m}[r]{$l_3$}
\epsfig{file=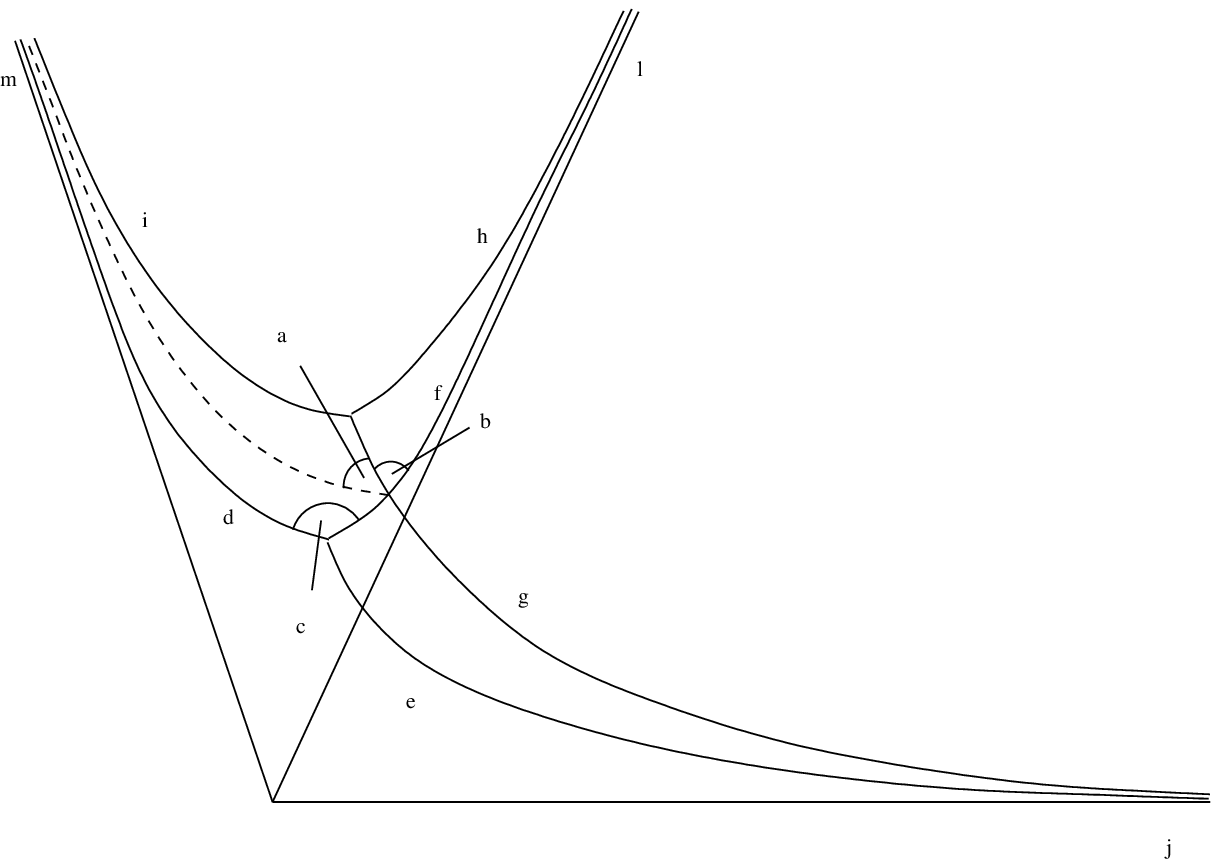, 
  width=0.8\textwidth}
\caption{}
\end{figure} 

\noindent We can use this monotonicity to prove uniqueness of a
homothetically expanding network. Let
$l_1, l_2, l_3$ be three given distinct half-lines, meeting at the origin and assume that there
exist two different homothetically expanding networks, satisfying the
balancing condition at the triple point such that both are asymptotic to the
three half-lines. Let the first network consist of the three curves
$b_1,b_2,b_3$, meeting at the triple point $p$ and 
the second consist of $c_1, c_2, c_3$, meeting at the triple point $q$. To avoid
notational confusion, let us assume that we have a configuration as in
Figure 2. It is easily seen that this argument works for any
such configuration. In the case of a degenerate configuration,
i.e.\;for example if $\text{Im}(b_1)\subset \text{Im}(c_1)$, it follows
easily from Lemma \ref{angle mon} that we get a contradiction along $c_1$.
In the picture, we have introduced the additional self-similarly expanding
curve, asymptotic to $l_3$ (dotted in the picture), starting at the intersection point $r$ of
$b_2$ and $c_1$. This curve intersects the
curve $c_1$ at an angle $\eta$ as shown above. The angle
between $b_2$ and $c_1$ at $r$ is denoted by $\gamma$. Note that by
Lemma \ref{angle mon}, applied along $c_1$, between $q$ and $r$, we
have that $\eta <\pi/3$, and $\gamma<\pi/3$. Thus $\eta+\gamma < 2
\pi/3$. Applying the monotonicity now along $b_2$, between $r$ and
$p$, we see that also $\zeta < 2\pi/3$ which yields a
contradiction. This completes the proof of Theorem \ref{straight thm}.

\begin{appendix}
\section{Further Results and Open Problems}\label{quadrupel sec}

\noindent Assume that we have again three distinct half-lines $l_1,l_2,l_3$,
meeting at the origin, and let $P(l_1,l_2,l_3)$ be the triple point of
the unique homothetically expanding network which is asymptotic
to the three half-lines as in Theorem \ref{straight thm}. We want
to argue that $P(l_1,l_2,l_3)$ depends continuously on
$l_1,l_2,l_3$. Note that the map $\Phi$ which maps the triple point $P$
and the angle $\vartheta$ of the three tangent directions in space 
at $P$ to the three
asymptotic half-lines is continuous in $P$ and $\vartheta$ by the
continuous dependence on initial data and the exponential convergence
of the self-similarly expanding curves to half-lines. Theorem \ref{straight
  thm} implies that this map is one-to-one and onto. By Lemma \ref{big
  slope}, for $|P|$ sufficiently large, at least two of the half-lines are
close to each other. Thus $\Phi$ is proper. So we can again apply \cite[Satz
8.12 and Satz 8.24]{Querenburg} to see that $\Phi$ is a homeomorphism. We
obtain:

\begin{lemma}\label{cont dep lem}
The triple point $P$ of a homothetically expanding network depends
continuously on the three distinct asymptotic half-lines $l_1,l_2,l_3$.
\end{lemma}

\begin{figure}[htb]
\psfrag{a}{$l_1$}   
\psfrag{b}{$l_3$}   
\psfrag{c}{$l_2$}   
\psfrag{d}{$l^\prime_2$}
\epsfig{file=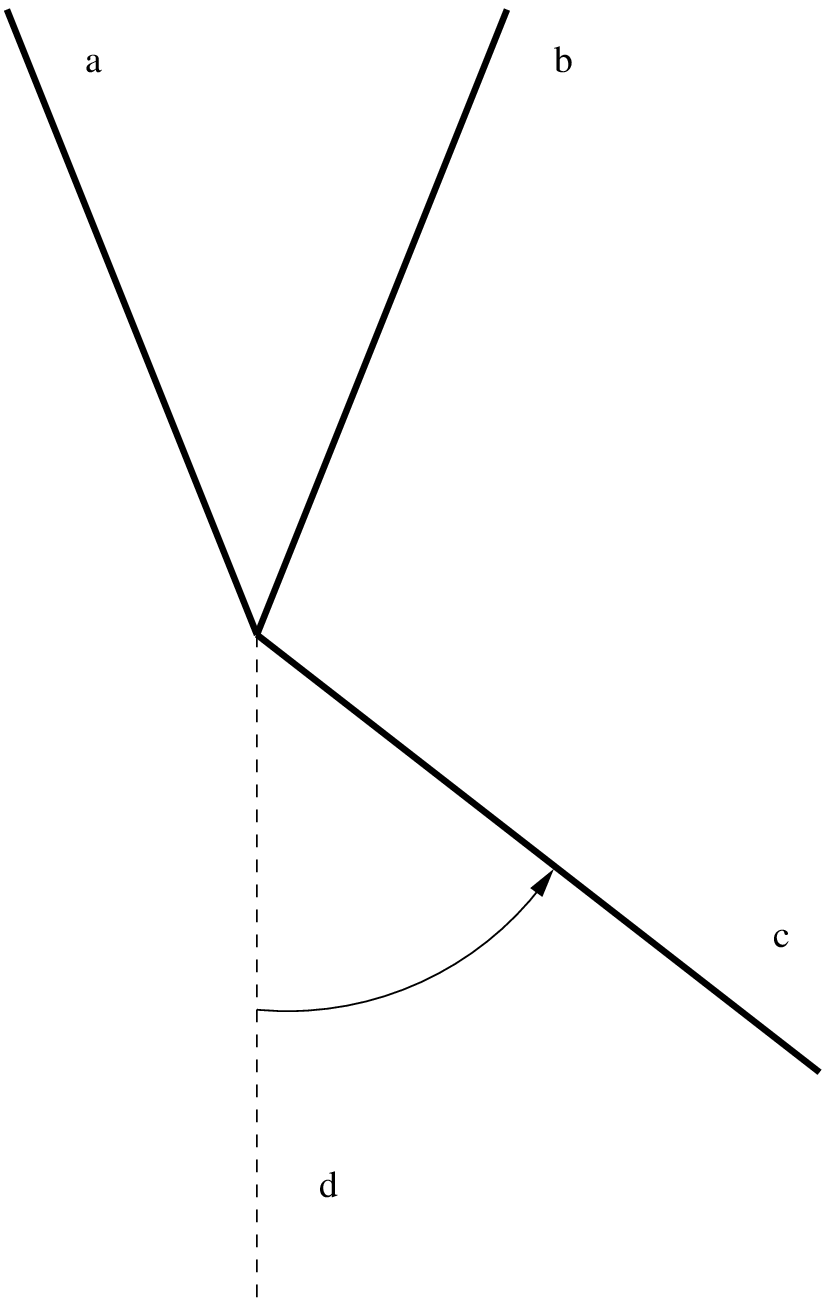, 
  width=0.3\textwidth} 
\caption{}
\end{figure}

\noindent Note that by uniqueness the triple point $P$ can lie on one of the
half-lines, say $l_2$, only if the configuration is symmetric, i.e.\! the
angle between $l_1$ and $l_2$ equals the angle between $l_2$ and $l_3$
and both are less or equal than the angle between $l_1$ and
$l_3$.
Furthermore, if the angles between $l_1$ and $l_2$ and between $l_2$ and $l_3$
are again equal, but bigger than the angle between $l_1$ and $l_3$,
then the triple point
is strictly contained in the smallest segment bounded by $l_1$
and $l_3$, see Lemma \ref{different slopes}. More precisely, it lies
on the continuation of $l_2$ into the segment bounded by $l_1$ and
$l_3$.

\noindent Now assume that we have a non-symmetric configuration
$l_1,l_2,l_3$, and the smallest segment is bounded by $l_1$ and
$l_3$. We want to show that the triple point has to lie in this
smallest segment, see Figure 3. To see this we first rotate $l_2$ into the symmetric
configuration, such that the segment bounded by $l_1$ and $l_3$
remains the smallest segment. Here we know that the triple point lies
strictly inside the segment bounded by $l_1$ and $l_3$. Then we
rotate $l_2$ back into its original position. Note that while doing
this we do not pass a symmetric configuration. Thus by Lemma \ref{cont
  dep lem} the triple point remains in the smallest segment, which gives:

\begin{lemma} Let $l_1,l_2,l_3$ be three distinct half-lines meeting
  at the origin, such that one of the segments bounded by the
  half-lines is the unique smallest one. Then $P(l_1,l_2,l_3)$ lies in
  this smallest segment.
\end{lemma}

\begin{remark} A homothetically expanding network as above partitions
  $\R^2$ into three unbounded segments. Note that locally the change
  of enclosed volume under mean curvature flow is given by the
  integral over the (mean) curvature along the curve. In this non-compact
  setting, taking for example the network at $t=0$ as a reference, the
  change of area of each of the sectors is given by the asymptotic
  angle of the sector minus $2\pi/3$.
\end{remark}

\begin{remark}
T. Ilmanen \cite{Ilmanen98} has studied the
existence of smooth solutions which approximate two intersecting
lines for small times and expand homothetically under curve shortening
flow. In general, these solutions are not unique. As indicated
in the picture,
\begin{figure}[htb]
\epsfig{file=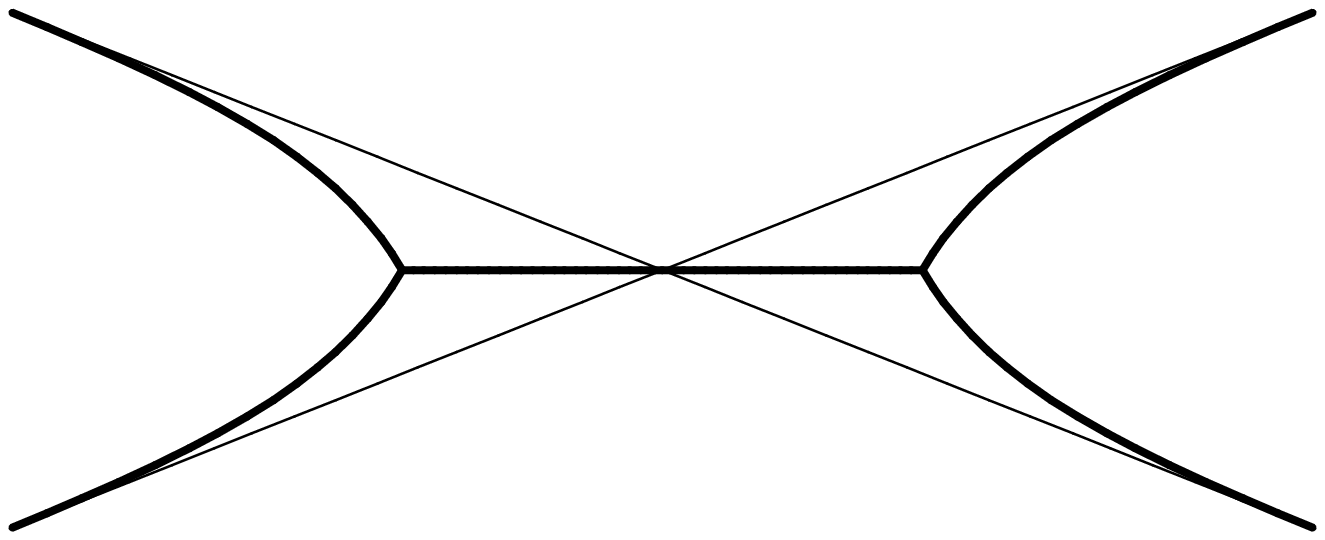, 
  width=0.4\textwidth}
\end{figure}
there exists another solution built from two networks as 
considered in this paper. The connecting segment is straight. 
Note that the triple points will be in sectors with an opening
angle of less than $120^\circ$. Such a construction can be
generalized, see also \cite{BrakkeBook}.
\end{remark}

\begin{problems}\rule{1ex}{0em}\\[-2ex]
\rule{0ex}{1em}
\begin{enumerate}[(i)]
\item Is there a simple (algebraic) relation involving the angles
between the half-lines and
\begin{itemize}
\item the direction, in which the triple point moves?
\item the tangent directions at the triple point?
\item the distance between the origin and the triple
point at $t=1/2$?
\end{itemize}
\item Show that there exists a solution to this initial value
problem if the network consists initially of three (not 
necessarily straight) curves with common start points.
\item Show that in this situation, the forward blow-up 
around the triple point is a solution as constructed in 
this paper. 
\end{enumerate}
\end{problems}
\end{appendix}

\bibliographystyle{amsplain}
\def\weg#1{}
\providecommand{\bysame}{\leavevmode\hbox to3em{\hrulefill}\thinspace}
\providecommand{\MR}{\relax\ifhmode\unskip\space\fi MR }
\providecommand{\MRhref}[2]{%
  \href{http://www.ams.org/mathscinet-getitem?mr=#1}{#2}
}
\providecommand{\href}[2]{#2}

\end{document}